\documentclass[10pt,a4paper,onecolumn]{article}
\usepackage{amsthm}
\usepackage{amsmath,amsfonts,amssymb,amsthm,epsfig,epstopdf,titling,url,array}

\usepackage{cite}
\theoremstyle{plain}
\newtheorem{thm}{Theorem}[section]

\newtheorem{prop}[thm]{Proposition}
\newtheorem{cor}{Corollary}

\theoremstyle{definition}
\newtheorem{defn}{Definition}[section]

\newtheorem{exmp}{Example}[section]
\newtheorem{rem}{Remark}[section]
\usepackage[english]{babel}
\usepackage[latin1]{inputenc}

\usepackage{amsfonts}
\usepackage{amssymb}
\usepackage{graphicx}
\usepackage{amsthm}
\begin{document}

	\author{
		Tareq Hamadneh$^{1}$,  Hassan Al-Zoubi$^{1}$,\\ Hamza Alzaareer$^{1}$, and Rafael Wisniewski$^{2}$ \\ \\
		Al Zaytoonah University of Jordan, Amman, Jordan$^{1}$\\
		Aalborg University, 9220 Aalborg East, Denmark$^{2}$\\ 
		\{t.hamadneh, dr.hassanz, h.alzaareer\}@zuj.edu.jo$^{1}$\\
		raf@es.aau.dk$^{2}$
		\\
	}
	\title{Optimization and Positivity Certificates of Rational Functions using Bernstein Form}
	\maketitle

	\begin{abstract}
		Rational functions of total degree $l$ in $n$ variables have a representation in the Bernstein form defined over $n$ dimensional simplex. The range of a rational function is bounded by the smallest and the largest rational Bernstein coefficients over a simplex. Convergence properties of the bounds to the range are reviewed. Algebraic identities certifying the positivity of a given rational function over a simplex are given. Subsequently, a bound established in this work does not depend on the given dimension.
		
	\end{abstract}
	
	
	\section{Introduction}
	The problem to decide whether a given rational function in $n-$variables is positive, in the sense that all its Bernstein coefficients are positive, goes back to [23, 24]. The same problem was addressed over different domains by other authors in [2] and [21]. The subject of certificates of positivity for polynomials over intervals was considered in [19, 20]. The same subject in the polynomial Bernstein basis was studied in [4], [14] and [21]. The problem of optimizing and approximating the minimum value of a polynomial over simplices was also extensively studied in [3], [5], [9, 13], [16]. The expansion of a given (multivariate) polynomial $p$ into Bernstein polynomials is used over a simplex, the so-called \textit{simplicial Bernstein form}, [3], [5], [7], [15, 16]. The approach was extended to the rational case in [8], [10] and [25], however, without investigating certificates of positivity for rational functions over simplices. Rational functions play an important role in stability analysis of dynamic systems, since the system can be stable with a rational control function [26]. Many researchers were focused on the topic of stability analysis for nonlinear systems, which continues to be a challenging problem [1, 18]. Recently, stability of nonlinear polynomial systems has been translated to certificates of positivity, [11, 12]. In this paper, we extend certificates of positivity to the multivariate rational Bernstein functions over simplices. This ensures the stability of nonlinear dynamic systems with a control in the rational Bernstein form [26]. Specifically, we will optimize the minimum and maximum values of the rational Bernstein function. Subsequently, we investigate certificates of positivity for rational functions in the Bernstein basis with respect to degree elevation and with respect to the maximum diameter of subsimplices. To this end, we will numerically provide valid bounds that approximate the number of subdivision steps and the degree of a rational function. Furthermore, we will provide a new technique for minimizing the range of a rational function over subsimplices. Finally, global and local certificates of positivity are satisfied with these bounds for rational functions in the Bernstein form.
	
	The organization of our paper is as follows. In the next section, we recall the most important background of the simplicial Bernstein expansion. In Section 3, we present polynomials in the Bernstein form. In Section 4, we extend optimization results to the rational case. Rational certificates of positivity are given in Section 5. Finally, Section 6 comprises conclusions.

	\section{Bernstein Expansion}
	We introduce some notation and necessary material about the simplicial Bernstein basis. Let $v_{0}, \dots,v_{n}$ be $n+1$ points of $\mathbb{R}^{n}  \;(n\geq 1)$, the ordered list $V=[v_{0}, \dots,v_{n}]$ is called simplex of vertices $v_{0}, \dots,v_{n}$. Throughout the paper, $V=[v_{0}, \dots,v_{n}]$ will denote a non-degenerate simplex of $\mathbb{R}^{n}$; viz the points $v_{0}, \dots,v_{n}$ are affinely independent. Let $\lambda_{0}, \dots, \lambda_{n}$ be the associated barycentric coordinates of $V$, i.e., the linear polynomials of $\mathbb{R}[X] = \mathbb{R} [X_{1}, \dots,X_{n}]$ such that $ \sum_{i=0}^{n} \lambda_{i}(x)=1$, and $ \forall x \in \mathbb{R}^{n}, \; x= \lambda_{0}(x) v_{0} + \dots+ \lambda_{n}(x) v_{n} $. The realization $|V|$ of the simplex $V$ is the subset of $\mathbb{R}^{n}$ defined as the convex
	hull of the points $v_{0},...,v_{n}$.
	
	We refer to the multi-index $ \alpha = ( \alpha_{0}, \dots,\alpha_{n} )  \in \mathbb{N}^{n+1}  $ and $ \lvert  \alpha    \lvert = \alpha_{0} + \dots + \alpha_{n}  $. Without loss of generality, we will often consider the standard simplex $ \Delta = [ e_{0}, e_{1}, \dots, e_{n}   ]  ,$ where $(e_{1}, \dots, e_{n} )$ denotes the canonical basis of $\mathbb{R}^{n}$, and $e_{0}= (0,...,0)$ the origin. This is not a restriction since any simplex $V$ in $\mathbb{R}^{n}$ can be mapped affinely upon $\Delta$. Subsequently, if $x =( x_{1}, ... , x_{n}) \in \Delta, $ then $(\lambda_{0}, ..., \lambda_{n})=  (  1 - \sum_{i=1}^{n} x_{i}, x_{1}, ..., x_{n}  )$. For $\hat{\beta}, \hat{\alpha} \in \mathbb{N}^{n}$ with $\hat{\beta} \leq \hat{\alpha}$, we define 
	\[  \binom{\hat{\alpha}}{\hat{\beta}} := \prod_{i=1}^{n} \binom{\alpha_{i}}{\beta_{i}}. \]
	
	If $k$ is any natural number such that $|\hat{\beta}|  \leq k$, we use the notation \[\binom{k}{\hat{\beta}}:= \frac{k!}{\beta_{1}!...\beta_{n}! (k-|\hat{\beta}|)! }.\]
	
	The Bernstein polynomials of degree $ k$ with respect to $V$ are the polynomials $(B_{\alpha}^{(k)})_{|\alpha|=k}, $ where
	
	\[
	B_{\alpha}^{(k)}(\lambda) = \binom{k}{\alpha} \lambda^{\alpha} .
	\]

	For $x \in \mathbb{R}^{n}$ its multipowers are $x^{\hat{\beta}}:= \prod_{i=1}^{n} x_{i}^{\beta_{i}} $.
	Let $p$ be a polynomial of degree $l,$
	
	\begin{equation}\label{egupoly}
	p(x) =  \sum_{|\hat{\beta}| \leq l}  a_{\hat{\beta}} x^{\hat{\beta}} ,
	\end{equation}
	$p$ can be uniquely expresses for $l\leq k$ as
	\[ p(x ) = \sum_{|\alpha| = k} b_{\alpha} (p,k,V) B_{\alpha}^{(k)} ,\]
	where $b_{\alpha} (p,k,V)$ are called the Bernstein coefficients of $p$ of degree $ k $ with respect to $V$.\\

	We recall the following notations:\\
	- The grid points of degree $k$ associated to $V$ are the points
	\[
	v_{\alpha}(k,V) = \frac{\alpha_{0} v_{0} + \dots + \alpha_{n} v_{n}  }{k} \in \mathbb{R}^{n} \; ( |\alpha| = k ) 
	\]
	which leads us to the control points associated to $p$ of degree $k$ with respect to $V$
	\begin{equation}
	( v_{\alpha}(k,V),\ b_{\alpha} (p,k,V)    )     \in \; \mathbb{R}^{n+1} \; ( |\alpha| = k   )    . 
	\end{equation}
	
	The control points of $p$ form its control net of degree $k$.\\
	- The discrete graph of $p$ of degree $k$ with respect to $V$ is formed by the
	points 
	\[  ( v_{\alpha}(k,V),\ p(v_{\alpha} (k,V)) )_{|\alpha| = k }       . \]
	
	\begin{prop}\label{keyfdrs}
		[16, Proposition 2.7] 
		For $p \in \mathbb{R}_{l} [X]$ and $k \geq l$, the following properties hold.\\
		
		\textbf{(i) Linear precision:} If degree $p \leq 1$, then
		\[ b_{\alpha} (p,k,V) = p (v_{\alpha} (k,V) ) , \ \forall  |\alpha| =k ; \]
		
		\textbf{(ii) Interpolation at the vertices:} If $(e_{0},...,e_{n})$ denotes the canonical
		
		basis of $\mathbb{R}^{n+1},$ then
		\begin{equation}\label{fp4}
		b_{ke_{i}} = p(v_{i}), \ 0 \leq i \leq n; 
		\end{equation}

		\textbf{(iii) Convex hull property:} The graph of $p$ over $V$ is contained in 
		
		the convex hull of its associated control points;\\
		
		\textbf{(iv) Range enclosing property:}
		\begin{equation}\label{fp3}
		\min_{|\alpha| = k} b_{\alpha} (p,k,V) \leq p(x) \leq \max_{|\alpha| = k}  b_{\alpha} (p,k,V) , \ \forall x \in V. 
		\end{equation}
	\end{prop}
	It follows from $(iv)$ in Proposition \ref{keyfdrs},
	the interval
	\[  B (p,k,V)   : =  [ \min b_{\alpha} (p,k,V) , \max b_{\alpha} (p,k,V)]  \]
	encloses the range of $p$ of degree $l \leq k$ over $V$.\\

	Finally, we denote the distance $d$ between two
	intervals $A=[\underline{a}, \overline{a}]$, $B=[\underline{b}, \overline{b}]$ by
	\[ d( A,B  )   := \max \{  |\underline{a}- \underline{b}|, |\overline{a}-\overline{b}|    \}  .   \]

	\section{Polynomial Bernstein Form}
	
	In this section, we present the most important properties of
	the Bernstein expansion over a simplex we will employ throughout the paper.
	

	In the following remark, we provide the simplicial polynomial Bernstein form of a given $p$ on $\Delta$. 
	\begin{rem}\label{p1}
		For $ \hat{\alpha}, \hat{\beta} \in \mathbb{N}^{n}$, let $p$ be a polynomial of degree $l$. The simplicial Bernstein form of $p$ of degree $l \leq k$  on $\Delta$ is given by
		\begin{equation}\label{fp1}
		p (x) =  \sum_{|\hat{\alpha}|+\alpha_{0} = k} b_{(\hat{\alpha},\alpha_{0})} (p,k,\Delta) B^{(k)}_{(\hat{\alpha},\alpha_{0})} (x) ,
		\end{equation}
		where
		\[	B^{(k)}_{ (\hat{\alpha}, \alpha_{0} ) } (x) =  \binom{k}{\hat{\alpha},\alpha_{0}}    x^{\hat{\alpha}} (1 - |x|  )^{\alpha_{0}}  ,   \ |\hat{\alpha}|+ \alpha_{0}  =k
		\]
		and
		\begin{equation}\label{fp2}
		b_{(\hat{\alpha},\alpha_{0})} (p,k,\Delta) = \sum_{\hat{\beta} \leq \hat{\alpha}  }   \frac{ \binom{\hat{\alpha}}{\hat{\beta}}   } {\binom{k}{\hat{\beta}}} a_{\hat{\beta}}      .   
		\end{equation}
	\end{rem}


	The Bernstein coefficients of degree $k$ ($l \leq k$) can be given as linear combinations of Bernstein coefficients of degree $l$, see, e.g., [15, Proposition 1.12].

	Let $(\hat{e}_{0},...,\hat{e}_{n})$ be points of $\mathbb{R}^{n+1}$, $\hat{e}_{i} = ( \underbrace{0,...,0}_{i},1, \underbrace{0,...,0}_{n-i}), \  i = 0,...,n$. By multiplying both sides of (\ref{fp1}) with $1 = (|x|+1 - |x|)^{k+1}$ and rearranging the result we obtain, see [6, Lemma 1.1] 
	\[ p =  \sum_{ |\beta |= k+1 } b_{\beta} (p,k+1,V) B_{\beta}^{k+1} , \]
	where
	\[   b_{\beta} (p,k+1,V) = \frac{1}{k+1} \sum_{i=0}^{n} 
	\beta_{i}  b_{\beta - \hat{e}_{i} } (p,k,V). \] 
	Hence, the range of $p$ of degree $k+1$ over $V$ can be bounded by
	\begin{equation}\label{fp0}
	B (p,k+1,V)   \subseteq B (p,k,V).
	\end{equation}
	
	\begin{rem}
		The number of Bernstein coefficients of an $n-$variate polynomial of degree $k$ is equal $\binom{k+n}{k}$.	
	\end{rem}
	
	The following definiton is given in [15].
	
	\begin{defn}
		Let $V = [ v_{0} , ..., v_{n} ]$ be a non-degenerate simplex of $\mathbb{R}^{n}.$ For $\gamma = k - 2$ and $0 \leq i < j \leq n,$ define the second differences of $p$ of degree $k$ with respect to $V$ as
		\[  \bigtriangledown ^{2}  b_{\gamma, i, j} (p,k,V) :=   b_{\gamma + e_{i} + e_{j - 1 }}   + b_{\gamma + e_{i-1} + e_{j}}  - b_{\gamma + e_{i-1} + e_{j - 1 }}  - b_{\gamma + e_{i} + e_{j  }} ,      \]
		with the convention $e_{-1}  := e_{n} . $ The second differences constitute the collection
		\[ \bigtriangledown^{2} b_{\gamma, i, j} (p,k,V) := ( \bigtriangledown^{2} b_{\gamma, i, j} (p,k,V)) _{ |\gamma| = k -2, 0 \leq i  < j\leq n .  }     \]
		
		Let $ || \bigtriangledown^{2} p ||_{\infty}$ denotes the maximum of the second differences, i.e.,
		\[ || \bigtriangledown^{2} p ||_{\infty} := \max_{ |\gamma|= k-2 , 0 \leq i  < j \leq n }  | \bigtriangledown^{2} b_{\gamma, i, j} (p,k,V)    | .  \]
	\end{defn}

	\begin{thm}\label{tp3}
		[15, Theorem 4.2] Let $ p \in \mathbb{R}_{l}[X]$ and $l <k$. Then
		\begin{equation}
		\max _{ |\alpha| = k  }  | p(v_{\alpha} (k,\Delta)  )  - b_{\alpha} (p,k,\Delta)  |  \leq  \frac{T}{k-1},
		\end{equation}
		where 
		\[ T := \frac{n(n+2)l (l-1)}{24} || \bigtriangledown^{2} p ||_{\infty} . \]
		
		A similar statement holds for the minimum. 	
	\end{thm}


	\section{Rational Bernstein Form}
	We may assume a rational function $f:=p/q$ where both $p$ and $q $ have the same degree $l$ since otherwise we can elevate the degree of the Bernstein expansion of either polynomial by component where necessary to ensure that their Bernstein coefficients are of the same order $ l \leq k.$ Since any simplex can be mapped upon the standard simplex by an affine transformation, we will extend results from polynomials to rational functions over $\Delta$.
	Let the range of $f$ over $\Delta$ be $ f(\Delta):= [\min f(x), \max f(x)] =:[ \underline{f}, \overline{f}]. $ The \textit{simplicial rational Bernstein coefficients} of $f$ of degree $k$ with respect to $\Delta$ are given by
	\begin{equation}\label{rnf1}
	b_{\alpha} (f,k,\Delta) = \frac{b_{\alpha} (p,k,\Delta)}{b_{\alpha} (q,k,\Delta)}, \ |\alpha|= k.
	\end{equation}
	Without loss of generality, we assume throughout the paper that  $b_{\alpha} (q,k,\Delta)>0$, $\forall |\alpha|=k$.

	The \textit{range enclosing property} for the rational function is given from [17, Theorem 3.1] as
	\begin{equation}\label{fp5}
	m^{(k)} := \min_{|\alpha| = k} b_{\alpha} (f,k,\Delta) \leq f(x) \leq \max_{|\alpha| = k}  b_{\alpha} (f,k,\Delta) =: M^{(k)} . 
	\end{equation}
	
	By application of (\ref{fp2}) to (\ref{rnf1}), the following theorem provides the sharpness property [25, Theorem 4] of $f$ with respect to its enclosure bound.
	
	\begin{thm}\label{ths1} 
		The equality holds in the right hand side of (\ref{fp5})
		\[
		\max_{x \in \Delta} f(x) =  M^{(k)}
		\]
		if and only if
		\[
		M^{(k)}  =    b_{(\hat{\alpha}^{*},\alpha_{0}^{*})} (f,k, \Delta) \ \text{for some} \ \hat{\alpha}^{*} = ke_{i_{0}} , \ i_{0} \in \{ 0,...,n  \} ,
		\]
		and $\alpha_{0}^{*} = k - |\hat{\alpha}^{*}|$.
		A similar statement holds for the equality in the left hand side of (\ref{fp5}).\end{thm}
	
	\begin{rem}
		We conclude from [25, Theorem 5] that
		\begin{equation}\label{frmn}
		B (f,k,\Delta)   \subseteq B (f,l,\Delta).
		\end{equation}
	\end{rem}

	In the following theorem, we review the linear convergence [25] of the range of a rational function to the enclosure bound under degree elevation. We include the positive and negative cases of $M^{(k)}$, since in [25] just the positive case is given.
	
	\begin{thm}\label{t1}
		For $l < k $ it holds that
		\begin{equation}
		d(f(\Delta) , B(f,k,\Delta) ) \leq \frac{\omega}{k-1} , 
		\end{equation}
		where	
		\begin{equation}\label{pfn}
		\omega :=      \frac{n  (n+2) l  (l-1)   }{24 \min_{ |\alpha| = l  }   b_{\alpha} (q,l,\Delta) } ( ||  \bigtriangledown^{2} p ||_{\infty}     +  \zeta  ||  \bigtriangledown^{2} q ||_{\infty}  )  , 
		\end{equation}
		and
		\begin{equation}\label{6fn27}
		\zeta  = \max \{ |\min_{ |\alpha| = l  }  b_{\alpha} (f,l,\Delta)|, |\max_{ |\alpha| = l  }  b_{\alpha} (f,l,\Delta)| \}.
		\end{equation}
		
		
	\end{thm}

	\textit{Proof.} The proof follows by using arguments similar to that given in the proof of Theorem \ref{thm2}.\\

	Assume that $ \Delta $ has been subdivided with respect to a point $\hat{v}$ in $\Delta$, $ \Delta = V^{[1]}  \cup ... \cup V^{[\sigma]}$, where the interiors of the simplices $|V^{[i]}| \ (1 \leq i \leq \sigma)$ are disjoint. Denote the union of the enclosure bounds over $V^{[i]}$, $i=1,...,\sigma$, by $B (p,k,V^{[\Delta]})$.
	The following theorm reviews the quadratic convergence [25, Theorem 7] of the range of a rational function to the enclosure bound with respect to subdivision. The proof also follows by using arguments similar to that given in the proof of Theorem \ref{thm2}.
	\begin{thm} \label{ptn1} Let $\Delta = V^{[1]} \cup ...\cup V^{[\sigma]} $ be a subdivision of the standard simplex $\boldsymbol{\bigtriangleup}$ and $h$ be an upper bound on the diameters of the $V^{[i],}$s. Then we have for $i=1,...,\sigma$
		\begin{equation}
		d \big(f(\Delta), B(f,k,V^{[\Delta]}) \big)   \leq h^{2} \omega^{\prime},
		\end{equation}
		where
		\begin{equation}\label{fp99}
		\omega^{\prime}  :=   k \frac{n^{2} (n+1) (n+2)^{2}   (n+3)   }{ 576 \min b_{\alpha} (q,l,\Delta) } \big( ||  \bigtriangledown^{2} p ||_{\infty}     +  \zeta  ||  \bigtriangledown^{2} q ||_{\infty}  \big) ,
		\end{equation}
		and $\zeta$ is the constant (\ref{6fn27}) independent of $h$.
	\end{thm}

	\section{Rational Certificates of Positivity}
	We study the positivity of rational functions over a non-degenerate simplex $V$. In order to do so, we use the simplicial rational Bernstein form. Certifying the positivity of rational functions is desired in many applications such stability analysis of dynamic systems, optimization and control theory. The enclosure property of $f$ shows that if all Bernstein coefficients of $f$ over $V$ are positive, then the rational function $f$ is positive over $V$. The converse need not to be true. There are rational functions which are positive over $V$ and some Bernstein coefficients $b_{\alpha} (f,l,V)$ are negative.
	\begin{exmp}\label{exo}
		Let 
		\[f(x) =  \frac{7x^{2} -5x +1 }{ x^{2} -2x +7 }, \] which is positive over $[-1,1]$, but the list of Bernstein coefficients $b_{\alpha} (f,2,[-1,1]) = (1.3,-1,0.5)$ has a negative value at $b_{1}(f)$.
	\end{exmp}
	
	The (univariate) Bernstein polynomials of $p$ of degree $k$ on $[\underline{x}, \overline{x}]$ 
	
	\[ B_{i}^{(k)}(x)= \binom{k}{i } \frac{(\overline{x}- x )^{k-i} ( x - \underline{x} )^{i}  }{w(X)^{k}}, \ i = 0,...,k, \]   
	take positive values over $(\underline{x}, \overline{x})$. Note that $B_{0}^{(k)}(x)$ is positive at $\underline{x}$ and $B_{k}^{(k)}(x)$ is positive at $\overline{x}$. The Bernstein coefficient $b^{(k)}_{0} $ is the value of $p$ at $\underline{x}$ and $b^{(k)}_{k} $ is the value at $\overline{x}$.
	Hence, if all Bernstein coefficients of $f=p/q$ are positive, the rational Bernstein form of $f$ over a given domain provides certificates of positivity for $f$ over the same domain. Without loss of generality, we assume that the (multivariate) rational case is studied on the standard simplex $\Delta$. Denote by $b(f,l,\Delta)$ the list of Bernstein coefficients of a rational function $f$ with respect to $\Delta$, we define $Cert (b(f,l,\Delta))$ by:
	\[ Cert (b(f,k,\Delta)) :  \left\{ \begin{array}{ll}
	b_{\alpha}(f,k,\Delta) \geq 0   & \mbox{ $ \text{for all} \ |\alpha| = k$}\\
	b_{k_{e_{i}}}(f,k,\Delta) > 0   & \mbox{ $ \text{for all} \ i \in \{ 0,...,n  \} $},\end{array} \right.     \]
	The rational Bernstein form of $f$ of degree $k$ is positive on $\Delta$ if  $\min_{|\alpha| = k  }  b_{\alpha} (f,k, \Delta) > 0.  $ In the following subsections, we decide if a rational function is positive and gives certificates of positivity in the rational Bernstein form by sharpness, degree elevation (global certificates), subdivision (local certificates) and minimization of a rational function. At the last, we provide a bound does not depend on the number of variables of $f$. 
	
	\subsection{Certificates by Sharpness}
	The sharpness property in Theorem \ref{ths1} satisfies the certificate of positivity of a rational function over $\Delta$. The equality holds in the left hand side of (\ref{fp5}) if \[ \min_{ |\hat{\alpha}|+\alpha_{0} = k  }   b_{\hat{\alpha}} (f,k,\Delta)   =    b_{ke_{i}} (f,k,\Delta) \ \text{for some}  \ i \in \{ 0,...,n  \}. \]
	This implies the following proposition. 
	
	\begin{prop}
		Given $f $ is positive on $\Delta$. If $\min_{ |\hat{\alpha}|+\alpha_{0} = k  }   b_{(\hat{\alpha},\alpha_{0})} (f,k,\Delta)   =    b_{ke_{i}} (f,k,\Delta) \ \text{for some} \ i , \ i \in \{ 0,...,n  \}$, then $f$ satisfies the certificate of positivity.
	\end{prop}

	\subsection{Global Certificates}
	If $k \geq l$ big enough, the minimum rational Bernstein coefficient of $f$ converges linearly to the minimum range $\underline{f}$ over $\Delta$. We show that
	the positive rational function has a global certificate of positivity at degree $k$ over $\Delta$.  The Bernstein degree is estimated in the following theorem.
	\begin{prop}\label{pcn1}
		Given $f$ is a positive rational function of degree $l$ over $\Delta$.
		If
		\[ k > \frac{\omega}{\underline{f}} +1,  \]
		where $\omega$ is the constant (\ref{pfn}), then $f$ satisfies the global certificate of positivity.
		
	\end{prop}

	\textit{Proof.} Let $k \geq l$ so that
	\[ \underline{f} -  m^{(k)}  \leq \underline{f} .  \]
	Then $b_{\alpha} (f,k,\Delta)$ are nonnegative.
	Theorem \ref{t1} implies that
	\[  \underline{f} -m^{(k)} \leq \frac{\omega}{k-1}, \]
	the interpolation property shows that $b_{ke_{i}}$, $\forall i \in \{0,...,n\}$, are positive. $\ \Box$\\
	
	Observing the obtained global certificate of positivity, we give the following corollary.
	\begin{cor}\label{thmn21}
		If $f$ is a rational function of degree $l$ positive over $\Delta$, then there exist some $k \geq l$ such that the minimum rational Bernstein coefficient of $f$ of degree $k$ is positive.  
	\end{cor}
	\begin{exmp}
		Let a rational function
		\begin{equation}\label{lne}
		f(x) = \frac{5x^{2} - 3 x  +1}{x^{2}+1 }
		\end{equation}
		of degree $l=2$, which is positive over $[0,1]$. Note that $\min b_{\alpha}(f,2,[0,1]) = -0.5 $ is negative. The rational Bernstein form of $f(x)$ (\ref{lne}) has a global certificate of positivity at $k =3$, since $\min b_{\alpha}(f,3,[0,1]) = 0 $, $b_{0}(f,3,[0,1])=1$ and $b_{3}(f,3,[0,1])= 1.5$.
	\end{exmp}
	


	\subsection{Local Certificates}
	In this section, we will not elevate the degree any more. This will lead to local certificates of positivity. 
	
	\begin{defn} [15, Definition 5.4]
		Let $S(\Delta) = (V^{[1]},...,V^{[\sigma]})$ be a subdivision of the simplex $\Delta$, $i.e., \ \Delta = V^{[1]} \cup...\cup V^{[\sigma]} $ and the interiors of the simplices $|V^{[i]}|$ are disjoint. If $f$ satisfies the certificate of positivity $Cert(b (f,l,V^{[i]}))$ for all $i = 1 ,...,\sigma$, we say that $f$ satisfies the local certificate of positivity associated to the subdivision $S(\Delta)$, which we write $Cert (b (f,l,S(\Delta))). $   
	\end{defn}
	
	We recall that the subdivision scheme consisting in $\frac{n(n+1)}{2}$ steps of binary splitting has a shrinking factor $1/2$.
	
	\begin{rem}
		From [15, Definition 5.5] and [16, Definition 2.14], the mesh of $\Delta$, denoted by $\hat{m},$ is its diameter. If $S$ is a subdivision scheme, we write $S^{N}(\Delta)$ the subdivision of $\Delta$ obtained after $N$ successive subdivision steps. $S$ is said to have a shrinking factor $0  < C <1$ if for every simplex $\Delta$, $\hat{m}(S(\Delta)) \leq C \times \hat{m}(\Delta)$, where $\hat{m}(S(\Delta))$ is the largest mesh among the subsimplices $U^{[i]}$. 	
	\end{rem}
	From Theorem \ref{ptn1}, the following proposition can be similarly shown as proposition \ref{pcn1}.
	
	\begin{prop}\label{thm1}
		Let $f$ be a rational function, positive over $\Delta$. Let $N$ be an integer and $S$ a subdivision scheme with a shrinking factor $ C < 1$. 
		Assume that
		\[  \frac{1}{C^{N}} >   \frac{\sqrt{2 \omega^{\prime}}}{\sqrt{\underline{f}}} ,  \]
		where $\omega^{\prime}$ is the constant (\ref{fp99}). Then $f$ satisfies the local certificate of positivity associated to $S^{N}(\Delta)$.

	\end{prop}

	\begin{exmp}
		We consider the rational function $f$ (\ref{lne}) over $I=[-1,1]$. The coefficients of $f$ over subintervals of width $1/2$ are given as follows:
		
		$b_{\alpha} (f,2, [-1,-1/2])= (1.3, 0.91, 0.63)$, \ $b_{\alpha} (f,2, [-1/2,0])= (0.63, 0.3, 0.14)$,
		$b_{\alpha} (f,2, [0, 1/2])= (0.14, -0.03, 0.04)$  and $b_{\alpha} (f,2, [ 1/2, 1])= (0.04, 0.12, 0.5)$. The rational Bernstein function still has a negative value over $[0, 1/2]$. Therefore, halving the interval $ [0, 1/2]$ and finding the coefficients over the new subintervals of $[0, 1/2],$
		\[  b_{\alpha} (f,2, [0, 1/4])= (0.14, 0.05, 0.02) \ \text{and} \ b_{\alpha} (f,2, [ 1/4, 1/2])= (0.02, 0, 0.4) ,\] satisfy the local certificate of positivity at the second subdivision step. 
	\end{exmp}
	
	\subsection{Optimization of Rational Functions}
	The enclosure bound leads to a lower bound of $f$ on $\Delta$. By repeatedly subdividing $\Delta$, the minimum of $f$ over $\Delta$ can then be approximated within any desired accuracy. Choosing $l= k$, the number of subdivision steps is bounded in Theorem \ref{thm2}. 

	\begin{rem}
		
		Let a rational Bernstein form of $f$ be on $\Delta$ and let the minimum Bernstein coefficient of $f$ of degree $l$ be $m$. The value $\delta$ is defined as	
		\[  \delta  = \min \big \{  f(v_{\alpha^{*}} ) , \ b_{le_{i}}(f), \ i = 0,...,n \big \},\] 
		where $m $ is attained at $ \alpha^{*}  $, $|\alpha^{*}|=l$.
		Then by (\ref{fp4}) and (\ref{fp5}) one can deduce that
		
		\begin{equation}
		m \leq \min_{x \in \Delta} f(x) \leq \delta.
		\end{equation}
		
	\end{rem}

	\begin{thm}\label{thm2}
		Let a rational function $f$, and $S$ a subdivision scheme with a shrinking factor $C<1 $. Let $\epsilon> 0$ and $N$ an integer satisfying
		\[ \frac{1}{C^{2N}}   >       \frac{ 2 \omega ^{\prime}}{\epsilon}    , \]
		where $\omega^{\prime}$ is the constant (\ref{fp99}), then
		\begin{equation}
		| \delta - m   | < \epsilon  .
		\end{equation}
	\end{thm}

	
	
	

	\textit{Proof.} Assume that 
	\begin{equation}\label{mfrasom}
	m = \frac{b_{\alpha^{*}} (p,l,\Delta)}{b_{\alpha^{*}} (q,l,\Delta)} , \ \text{for some} \ |\alpha^{*}|=l .
	\end{equation}
	Let $\zeta  = \max \{ |\min_{ |\alpha| = l  }  b_{\alpha} (f,l,\Delta)|, |\max_{ |\alpha| = l  }  b_{\alpha} (f,l,\Delta)| \}.$
	We can conclude from (\ref{mfrasom}) and the corresponding grid point $v_{\alpha^{*}}$ that
	\[  \delta - m  \leq  f(v_{\alpha^{*}}) - b_{\alpha^{*}} (f,l,\Delta)   \]
	\[ =   \frac{ p(v_{\alpha^{*}}) - m . q(v_{\alpha^{*}})   + m . b_{\alpha^{*}} (q,l,\Delta)   - b_{\alpha^{*}} (p,l,\Delta) }{q(v_{\alpha^{*}}) }     \]
	\[ \leq  \frac{ p(v_{\alpha^{*}}) - b_{\alpha^{*}} (p,l,\Delta) + m   . b_{\alpha^{*}} (q,l,\Delta) - q(v_{\alpha^{*}})       }{ \min b_{\alpha} (q,l,\Delta)}     .  \]
	Taking absolute values and using (\ref{fp5}) we can estimate
	\[ \delta - m  \leq C^{2N} l  \frac{n^{2} (n+1) (n+2)^{2} (n+3)  }{ 288 \min b_{\alpha} (q,l,\Delta)} \big(    ||  \bigtriangledown^{2} p ||_{\infty}     +  \zeta  ||  \bigtriangledown^{2} q ||_{\infty}         \big)    ,  \]
	where the last inequality follows by Proposition \ref{thm1} which completes the proof.
	
	\begin{cor} Given $f$ is a rational function of degree $l$ over $\Delta$. Assume under the assumptions of Theorem \ref{thm2} that $ \underline{f} \geq \epsilon $. If $N$ and $\underline{f} $ are satisfying:
		
		\[\frac{1}{C^{2N}} >  \frac{ 2 \omega^{\prime}}{\underline{f}} , \]
		then $f$ satisfies the local certificate of positivity $Cert (b (f,l,S(\Delta))). $ 
	\end{cor}

	\subsection{Independent Bounds}
	In this section, we provide a bound does not depend on the number of variables of $f=p/q$. Such this bound is the best in high dimensions as explained in [9], [14] and [22].

	Powers and Reznick in [22] have proved the following bound:
	\begin{thm} \label{thm6r3} [21, Theorem 3]
		Let $p$ be a polynomial of degree $l$, positive on the standard simplex $\Delta$. Let $\underline{p}$ be the minimum of $p$ on $\Delta$. Then for
		\begin{equation}\label{fern}
		k > \frac{   l (l-1)  }{2} \frac{\max |b_{\alpha}  (p,l,\Delta)| }{\underline{p} },
		\end{equation}
		the Bernstein form of $p$ of degree $k$ has positive coefficients.
	\end{thm}
	
	In the following corollary, we hold Theorem \ref{thm6r3}, [22, Proposition 4] and results from [21] to the rational case. 
	\begin{cor}
		Let $f=p/q$ be a rational function, positive on $\Delta$. If	
		\[   k > \frac{   l (l-1)  }{2} \frac{\max |b_{\alpha}  (p,l,\Delta)| }{\underline{p} },  \]
		then $f$ satisfies the global certificate of positivity.
	\end{cor}
	
	\textit{Proof.} Let $k$ be large enough and $\frac{p(\boldsymbol{x})}{q(\boldsymbol{x})} > 0=: a $, hence $h(\boldsymbol{x}) := p(\boldsymbol{x}) - a  q(\boldsymbol{x}) > 0$ allows from Theorem \ref{thm6r3} and [22, Proposition 4] the (global) certificate of positivity. Hence, If	
	\[
	k > \frac{l(l-1)}{2} \frac{\max |b_{\alpha} (h,l,\Delta)|  }{ \underline{p} },
	\]
	then $b_{\alpha} (p,l,\Delta) -  a.b_{\alpha} (q,l,\Delta)  > 0 $.
	It follows that $b_{\alpha} (h,l,\Delta)/b_{\alpha} (q,l,\Delta)$ are positive, $\forall |\alpha| = l,$ which completes the proof .

	
	\begin{cor} 
		Let $D_{1} = \frac{\omega}{ \min f(x)   }+1$ and $D_{2} = \frac{   l (l-1)  }{2} \frac{\max |b_{\alpha}  (p,l,\Delta)| }{\underline{p} } $, where $\underline{p}$ is the minimum of $p$ on $\Delta$. Then the positive rational function $f=p/q$ satisfies the (global) certificate of positivity over $\Delta$ if $k > \max \{ D_{1}, D_{2}  \}.$
	\end{cor}
	\begin{rem}
		Given $f$ is a rational function of degree $l$, negative on $\Delta$. Then $f$ satisfies certificates of negativity by applying the same arguments above to the upper bounds.
	\end{rem}
	
	\section{Conclusions}
	
	In this paper, we considered the multivariate rational function $f = p/q$ in the Bernstein form. The expansion of the numerator and denominator polynomials into Bernstein polynomials was applied. We reviewed properties such sharpness and monotonicity of bounds of a multivariate rational function. The linear and the quadratic convergence of the enclosure bound to the range of a rational function $f$ improved the bounds of $f$ over simplices. By repeatedly subdividing the simplex, the minimum of $f$ over a simplex was approximated within a desired accuracy. We addressed an optimization of a (multivariate) rational function and bounded the number of subdivision steps. Subsequently, we investigated certificates of positivity in the simplicial rational Bernstein form by sharpness, degree elevation and subdivision. At the last, we estimated the degree of the Bernstein expansion by a bound which is not depending on the given dimension.\\

	\textbf{Acknowledgments} \\
	The authors gratefully acknowledge support from Al-Zaytoonah University of Jordan under the grant number 2019-2018/585/G12. The first author would like to thank Professor Amjed Zraiqat for his careful reading the manuscript and the constructive comments.\\
	




	\textbf{References}
	\begin{description}
 \item 1. Juergen Ackerman.  Robust control.  London:  Springer-Verlag, 1993.
 \item 2.	Richard Askey. Certain rational functions whose power series have positive coefficients. ii. SIAM Journal on Mathematical Analysis, 5(1):53-57, 1974.
 \item 3.	Saugata  Basu,  Richard  Leroy,  and  Marie-Francoise  Roy.  A  bound  on  the minimum of a real positive polynomial over the standard simplex. arXiv preprint arXiv:0902.3304, 2009.
\item 4.	Fatima Boudaoud, Fabrizio Caruso, and Marie-Francoise Roy.  Certificates of positivity in the Bernstein basis. Discrete and Computational Geometry, 39(4):639-655, 2008.
 \item 5.	Etienne De Klerk, Dick Den Hertog, and G Elabwabi. On the complexity of optimization over the standard simplex. European Journal of Operational Research, 191(3):773-785, 2008.
 \item 6.	Gerald  E  Farin.    Triangular  Bernstein-Bezier  patches.    Computer  Aided Geometric Design, 3(2):83-127, 1986.
 \item 7.	Rida T Farouki. The Bernstein polynomial basis: a centennial retrospective. Computer Aided Geometric Design, 29(6):379-419, 2012.
 \item 8. 	Juergen Garloff and Tareq Hamadneh. Convergence and inclusion isotonicity of the tensorial rational Bernstein form. In Scientific Computing, Computer Arithmetic, and Validated Numerics, pages 171-179. Springer, 2015.
\item 9. 	Tareq Hamadneh. Bounding Polynomials and Rational Functions in the Tensorial and Simplicial Bernstein Forms. PhD thesis, University of konstanz, 2018.
 \item 10.  Tareq Hamadneh, Nikolaos Athanasopoulos, and Mohammed Ali. Minimization and positivity of the tensorial rational Bernstein form. In 2019 IEEE Jordan International Joint Conference on Electrical Engineering and Information Technology (JEEIT), pages 474-479. IEEE, 2019.
 \item 11.	Tareq Hamadneh and Rafael Wisniewski. Algorithm for Bernstein polynomial control design. IFAC-PapersOnLine, 51(16):283-289, 2018.
 \item 12.	Tareq Hamadneh and Rafael Wisniewski. The barycentric Bernstein form for control design. In 2018 Annual American Control Conference (ACC), pages 3738-3743. IEEE, 2018.
\item	13. Eldon Hansen. Global optimization using interval analysis the multi-dimensional case. Numerische Mathematik, 34(3):247-270, 1980.
\item	14. Richard  Leroy.  Certificats de positivite et minimisation polynomiale dans la base de Bernstein multivariee.  PhD thesis, Universite Rennes 1, 2008.
\item	15. Richard Leroy. Certificates of positivity in the simplicial Bernstein basis. working paper or preprint, 2009.
\item	16. Richard Leroy. Convergence under subdivision and complexity of polynomial minimization in the simplicial Bernstein basis. Reliable Computing, 17:11-21, Springer, 2012.
\item	17. Anthony Narkawicz, Juergen Garloff, Andrew P Smith, and Cesar A Munoz. Bounding the range of a rational functiom over a box. Reliable Computing, 17(2012):34-39, 2012.
\item	18. Sergey G Nersesov and Wassim M Haddad. On the stability and control of nonlinear dynamical systems via vector lyapunov functions. IEEE Transactions on Automatic Control, 51(2):203-215, 2006.
\item	19. Yurii Nesterov. Squared functional systems and optimization problems. In High performance optimization, pages 405-440. Springer, 2000.
\item	20. Victoria Powers and Bruce Reznick. Polynomials that are positive on an interval. Transactions of the American Mathematical Society, 352(10):4677- 4692, 2000.
\item	21. Victoria  Powers  and  Bruce  Reznick.   A  new  bound  for  polya's  theorem with applications to polynomials positive on polyhedra. Journal of Pure and Applied Algebra, 164(1):221-229, 2001.
\item	22. Victoria Powers and Bruce Reznick. Polynomials positive on unbounded rectangles. In Positive polynomials in control, pages 151-163. Springer, 2005.
\item	23. Armin Straub and Wadim Zudilin. Positivity of rational functions and their diagonals. Journal of Approximation Theory, 195:57-69, 2015.\\
\item	24. Gabor Szego.  Ueber gewisse potenzreihen mit lauter positiven koeffizienten. Mathematische Zeitschrift, 37(1):674-688, 1933.
\item	25. Jihad Titi, Tareq Hamadneh, and Juegen Garloff.  Convergence of the simplicial rational Bernstein form. In Modelling, Computation and Optimization in Information Systems and Management Sciences, pages 433-441. Springer, 2015.
\item	26. Mohsen Vatani and Morten Hovd. Lyapunov stability analysis and controller design for rational polynomial systems using sum of squares pro- gramming. In 2017 IEEE 56th Annual Conference on Decision and Control (CDC), pages 4266-4271. IEEE, 2017.

\end{description}

	\medskip

\end{document}